\documentclass[aps,pre,longbibliography,twocolumn,superscriptaddress,showkeys,floatfix]{revtex4-2}

\usepackage{dcolumn}% Align table columns on decimal point
\usepackage{float}  %Display figures in section by using the H modified for the figure float
\usepackage{amsmath,amssymb,amsfonts}
\usepackage{bm}
\usepackage{graphicx}
\usepackage{hyperref}
\usepackage{xcolor}
\usepackage{booktabs}
\hypersetup{  % Remove ugly boxed around links
	colorlinks,
	linkcolor={red!50!black},
	citecolor={blue!50!black},
	urlcolor={blue!80!black}
}

\newcommand{\eml}{\operatorname{eml}}
\newcommand{\dd}{\mathrm{d}}
\newcommand{\ee}{\mathrm{e}}

\begin{document}
	
	\title{Non-Monotone Response Modules and Cascades from the EML Operator for Reduced Models of Biological Dynamics}
	
	\author{Amir Erez}
	\email{amir.erez1@mail.huji.ac.il}
	\affiliation{Racah Institute of Physics, The Hebrew University of Jerusalem, Jerusalem 9190401, Israel}
	
	\date{\today}
	
	\begin{abstract}
		Standard saturating response functions, such as the Hill function, are monotone and therefore cannot represent recruitment-induced overshoot or adaptive transients with a single block. Reproducing such non-monotone responses from saturating primitives requires at least a difference of two blocks with opposing amplitudes, doubling the static-block parameter count. Here, building on a recent mathematical result that a single binary operator, EML, generates all standard elementary functions, we use EML as a structured grammar for reduced nonlinear ODEs. This yields an activation-suppression module that captures overshoot directly. We validate the framework in three settings. First, on PKA-R relocalization data, the EML grammar discovers a reduced surrogate consistent with established mechanistic biology. Second, on Rho-GTPase recruitment data, an exhaustive search over EML expression trees selects the same compositional form across all four perturbation-response traces. Third, a $50$-state simulated network is compressed by an EML cascade acting as a fixed temporal basis. Thus we demonstrate the power and potential of EML for reduced models of biological dynamics.
	\end{abstract}
	
	\maketitle
	
	\section{Introduction}
	
	Differential equations occupy a central position in many branches of physics and quantitative biology. They are useful to identify a small number of state variables, couplings, and time scales that organize the observed behavior \cite{transtrum2015,lopez2026}. In biological modeling, however, this reduction is often quite approximate. Perturbation experiments can reveal overshoot, adaptation, delayed inhibition, or biphasic dependence, while the underlying biochemical network may be only partially known. In such cases, one often falls back on standard response functions, such as Hill-type saturation \cite{nelson2022}, or on large mechanistic models whose parameters are difficult to identify. Here, we explore a controlled way to generate low-dimensional nonlinear ODEs, that is richer than a fixed phenomenological curve but more constrained than arbitrary curve fitting.
	
	Symbolic regression addresses a related problem: given data, search for an analytic expression that explains it. Modern approaches have shown that scientific laws and compact physical relations can sometimes be recovered from numerical data when the search is constrained by appropriate inductive biases \cite{schmidt2009,udrescu2020,cranmer2020}. Dynamical variants, including sparse identification of nonlinear dynamics, similarly seek compact governing equations from time-series data \cite{daniels2015,brunton2016,rudy2017}. Underlying these methods is a useful principle: model discovery is more meaningful when the search grammar, the model-complexity penalty, and the validation criteria are specified in advance. Otherwise, a sufficiently expressive symbolic system becomes a language for overfitting.
	
	Here, we explore a new grammar for reduced nonlinear dynamics based on the recently proposed EML operator by Odrzywo{\l}ek \cite{odrzywolek2026}, where the single binary operator
	
	\begin{equation}
		\eml(x,y)=\exp(x)-\ln(y)\,,
	\end{equation}
	
	\noindent together with the constant $1$ can generate the standard repertoire of elementary functions. In Odrzywo{\l}ek's work, the main point is syntactic universality: elementary formulae can be represented as binary trees with internal nodes from the EML operator. Thus, a natural grammar for symbolic regression emerges. Rather than searching over many disjoint primitives, e.g., power, exponential, logarithm, and arithmetic operations, one can search over EML-generated expression trees, easily translatable into ordinary mathematical form.
	
	We stress that EML is not a new physical or biological interaction law. A universal expression grammar supplies syntax, not semantics. Physical or biological interpretation enters only after the grammar has been restricted and the resulting model has been compared with simpler null models and, when available, mechanistic alternatives. In this manuscript, EML is used as a generator of candidate equations. The selected equations are then read as one would typically, with terms for activation, relaxation, and hidden coarse-grained states.	
	
	We define a simple yet useful EML-generated module,
	
	\begin{equation}
		\label{eq:EMLmodule}
		M_1(R) = \eml(\alpha\ln R,\ee^{\beta R})=R^\alpha-\beta R .
	\end{equation}
	
	This expression is a minimal candidate for adaptive or biphasic responses: increasing input can initially increase the output, while sufficiently strong input can suppress it. Written in expanded form, the expression is more readily interpretable than the EML notation itself. The role of EML is to provide the compositional rule by which this and higher-order modules are generated.
	
	The well-known Hill functions remain an appropriate null model for monotone saturating responses \cite{hill1910}. They are interpretable and widely used when modeling biological systems. The EML-generated module does not attempt to replace them when they work. Its purpose is to provide an alternative when the observed response is nonmonotone, adaptive, or shaped by competing positive and negative processes \cite{Vogel2016}. In such cases, a monotone response function is insufficient, whereas our EML module, Eq.~\ref{eq:EMLmodule}, supports an internal optimum with a small number of parameters. 
	
	Our central claim is methodological. The EML gate is non-monotone in its argument with three parameters, and is therefore the minimal activation-suppression block that can produce a rise-then-fall response under monotone input. Standard saturating primitives, e.g.\ the Hill function, are monotone and require at least a difference of two such blocks to produce the same shape. Building on this depth-one structural asymmetry, we use EML as a compact, differentiable, and recursively expandable grammar for constructing reduced nonlinear ODEs: at depth one it generates activation-suppression response functions; at higher depth it generates cascades of hidden states whose solutions are given by sequential convolution integrals. We derive a controlled hierarchy of models in which one can increase EML depth and ask whether the additional structure is justified by held-out predictive error. Importantly, EML depth is not a count of molecular steps; it is an effective coarse-grained dynamical depth, useful when hidden delays or adaptive processes shape the measured response.
	
	This manuscript develops our methodology in three steps. First, we apply an EML-ODE grammar to the PKA regulatory-subunit relocalization data of LaCroix \emph{et al.} \cite{lacroix2022}, who measured and then proposed a mechanistic linker-occupancy model for the observed time-series. The EML grammar is therefore tested as a surrogate, not as a replacement for the known mechanism. Second, we apply the restricted grammar to Rho-GTPase perturbation data from Nanda \emph{et al.} \cite{nanda2023} and compare it directly with a parallel Hill-grammar search of equal compositional depth, making the depth-one asymmetry between the two grammars concrete. Third, we construct a toy high-dimensional activation-adaptation network and ask whether a low-dimensional EML cascade can capture its behavior.

	\section{Results}
	
	We start from the activation-suppression module, Eq.~\ref{eq:EMLmodule}. Depending on the system, the suppression term may represent sequestration, adaptation, toxicity, or another process opposing the positive first term. For $0<\alpha<1$ and $\beta>0$, $M_1(R)$ is unimodal. The optimal input is $R_*=\left(\frac{\alpha}{\beta}\right)^{1/(1-\alpha)}$, derived in Appendix~\ref{app:module}.	Thus, a monotonically increasing input can produce a transient overshoot if it passes through the optimal point. This feature distinguishes the EML activation-suppression module from a monotone Hill response, which can rise and saturate or fall and saturate but cannot rise and then fall under a monotone input without adding a second term. The minimal time-dependent response model is a first-order relaxation equation,
	\begin{align}
		\label{eq:K1_ode}
		\tau \dot y
		&= -y+y_0+B M_1(R(t)), \nonumber \\
		&= -y+y_0+B\left[R(t)^\alpha-\beta R(t)\right].
	\end{align}
	For chemically induced recruitment experiments, a useful input model is
	\begin{equation}
		\label{eq:R_recruitment}
		R(t)=R_\infty\left(1-\ee^{-k_R t}\right),
	\end{equation}
	with $R_\infty$ proportional to the applied perturbation strength. Equations~\eqref{eq:K1_ode} and \eqref{eq:R_recruitment} define a minimal activation-suppression ODE. It has the convolution solution (Appendix~\ref{app:solutions}),
	\begin{equation}
		\label{eq:K1_solution}
		y(t)=y(0)\ee^{-t/\tau}
		+\frac{1}{\tau}\int_0^t \ee^{-(t-s)/\tau}
		\left[y_0+B M_1(R(s))\right]\,\dd s .
	\end{equation}

	\subsection{EML as a symbolic-regression grammar}
	\label{sec:symbolic_regression}

	The preceding equations can be read as an explicit reduced model, or, as the depth-one member of a restricted EML expression grammar which defines our proposed methodology. Below, we generate candidate models from $R(t)$, binary sums, and a centered EML gate
	\begin{equation}
		\label{eq:centered_gate}
		G_{a,b,c}(x)=\eml\!\left(a\ln[c+x],\ee^{bx}\right)-c^a
		=(c+x)^a-bx-c^a .
	\end{equation}
	The subtraction centers the gate at zero input and is not a new mechanism; it fixes the baseline of the generated expression. The restricted grammar is specified compactly as
	\begin{equation}
		\label{eq:grammar}
		E ::= R \mid G(E) \mid E+E .
	\end{equation}
	For the unfamiliar reader: the symbol $E$ denotes a candidate expression; ``$::=$'' should be read as ``is generated as''; vertical bars ``$\mid$'' mean ``or''. Thus, Eq.~\eqref{eq:grammar} says that a candidate expression $E$ may be one of three things: the input variable $R$ itself; an EML gate applied to a previously generated expression, $G(E)$; or the sum of two previously generated expressions, $E+E$. For example, the grammar first allows $R$. Applying the second rule gives $G(R)$. Applying it again gives $G(G(R))$. Using the addition rule gives expressions such as $G(R)+R$ or $G(R)+G(R)$. Applying a gate to a sum gives expressions such as $G(G(R)+R)$. These are the expression trees that are fitted and compared. To avoid duplicates, $G(R)+R$ and $R+G(R)$ are treated as the same expression. A fitted static response, i.e., an instantaneous algebraic mapping without intrinsic relaxation time, has the form
	\begin{equation}
		\label{eq:symbolic_response}
		y(t)=y_0+B\,E(R(t;k)),
		\qquad
		R(t;k)=1-\ee^{-kt}.
	\end{equation}
	The search then enumerates all expressions satisfying prescribed depth and node limits, fits their continuous parameters by nonlinear least squares on training time-series, and ranks models by held-out weighted error together with optional penalties on expression depth and node count. For observations $y_i$ with standard errors $\sigma_i$, training set $\mathcal T$, and held-out set $\mathcal V$, the fitted parameters are
	\begin{equation}
		\hat\theta_E=\arg\min_{\theta}
		\sum_{i\in\mathcal T}
		\left[
		\frac{y_i-\hat y_E(t_i;\theta)}{\tilde\sigma_i}
		\right]^2 ,
	\end{equation}
	where $\tilde\sigma_i=\max(\sigma_i,\sigma_{\rm floor})$ prevents zero or very small error values from dominating the fit. Model ranking uses
	\begin{equation}
		\label{eq:validation_score}
		\mathcal S(E)=
		\frac{1}{|\mathcal V|}
		\sum_{i\in\mathcal V}
		\left[
		\frac{y_i-\hat y_E(t_i;\hat\theta_E)}{\tilde\sigma_i}
		\right]^2
		+\lambda_d d(E)+\lambda_n n(E),
	\end{equation}
	where $d(E)$ is the EML depth and $n(E)$ is the number of expression-tree nodes. The fits shown here use $\lambda_d=\lambda_n=0$ and report the unpenalized held-out weighted mean squared error; AIC and BIC are also computed as diagnostic quantities but are not the primary selection criterion. Although the model comparison can be done differently, our methodological purpose is insensitive to a precise comparison regime: EML is not the explanation of a biological time-series. It is the engine used to propose compact equations. Therefore, the empirical examples below should be read as demonstrations of a constrained model-discovery framework rather than as claims that EML is a physical interaction law.

	A natural comparator is a Hill recruitment response
	\begin{equation}
		\label{eq:hill_static}
		H(R)=A\frac{R^h}{K_d^h+R^h},
	\end{equation}
	with the corresponding relaxation model
	\begin{equation}
		\label{eq:hill_ode}
		\tau \dot y=-y+y_0+H(R(t)).
	\end{equation}
	Equation~\eqref{eq:hill_ode} is the correct null model when the perturbation produces a monotone approach to a saturated activity. It is structurally inappropriate, however, for recruitment-induced overshoot or biphasic responses, because $H(R)$ is monotone in $R$ for every choice of $h>0$, $K_d>0$. Reproducing a rise-then-fall response from saturating primitives therefore requires at least a difference of two such blocks with opposing amplitudes, doubling the static parameter count. By contrast, the centered EML gate of Eq.~\eqref{eq:centered_gate} is non-monotone in $x$ for $0<a<1$, $b>0$ and produces an activation-suppression response with a single block. This single-block non-monotonicity is the structural asymmetry exploited in the experimental examples below, and we will return to it explicitly when we compare the EML grammar to a Hill grammar of equal compositional depth on the Nanda \emph{et al.} data (Sec.~\ref{sec:nanda_grammar}).
	
	\subsection{PKA-R relocalization in LaCroix et al.}
	\label{sec:lacroix_grammar_ode}

	The first example is the PKA regulatory-subunit relocalization experiment of LaCroix \emph{et al.} \cite{lacroix2022}, where rapamycin-induced recruitment of PKA-R to the plasma membrane produces a paradoxical response: low or intermediate PKA-R recruitment enhances plasma-membrane PKA activity, whereas high recruitment inhibits it. The authors interpret this as a linker or scaffold-like stoichiometric effect. At moderate abundance, PKA-R helps assemble productive local signaling complexes; at excessive abundance, it titrates cAMP or PKA-C away from productive complexes and suppresses activity.

	This dataset is useful because it has an independently motivated mechanism, where signaling is proportional to the concentration of fully occupied linker molecules \cite{lacroix2022}. A dynamical version of that model can be written as
	\begin{equation}
		\label{eq:lacroix_linker_ode}
		\tau \dot y_D
		=-y_D+1+A\left[\Phi_N(S_D(t))-\Phi_N(S_0)\right],
	\end{equation}
	where $S_D(t)=S_0+qD(1-\ee^{-k_Rt})$ and $\Phi_N(S)$ denotes the fully occupied linker concentration up to scale. For the $N=4$ case relevant to four cAMP-binding sites,
	\begin{equation}
		\label{eq:Phi4}
		\Phi_4(S)=\frac{S}{S^4+S^3+2S^2+3S+4}.
	\end{equation}
	The derivation and normalization of Eq.~\eqref{eq:Phi4} are summarized in Appendix~\ref{app:linker}.

	We used this example to test whether the restricted EML grammar can find a reduced surrogate. The search used the same expression grammar as Eq.~\eqref{eq:grammar}, but embedded the selected expression in a first-order kinetic equation rather than fitting it as an instantaneous static response:
	\begin{equation}
		\label{eq:lacroix_grammar_ode}
		\tau \dot y_D
		=
		-y_D+1+B\,E(R_D(t;k)),
		\qquad
		R_D(t;k)=D(1-\ee^{-kt}).
	\end{equation}
	This kinetic embedding is important. A purely static grammar, $y_D(t)=1+B E(R_D(t))$, can produce artificial early features because it forces PKA activity to follow recruitment instantaneously. The relaxation form in Eq.~\eqref{eq:lacroix_grammar_ode} removes that artifact and gives the grammar the same minimal kinetic structure as the Hill and linker comparators.

\begin{figure*}[t]
		\centering
		\includegraphics[width=1.98\columnwidth]{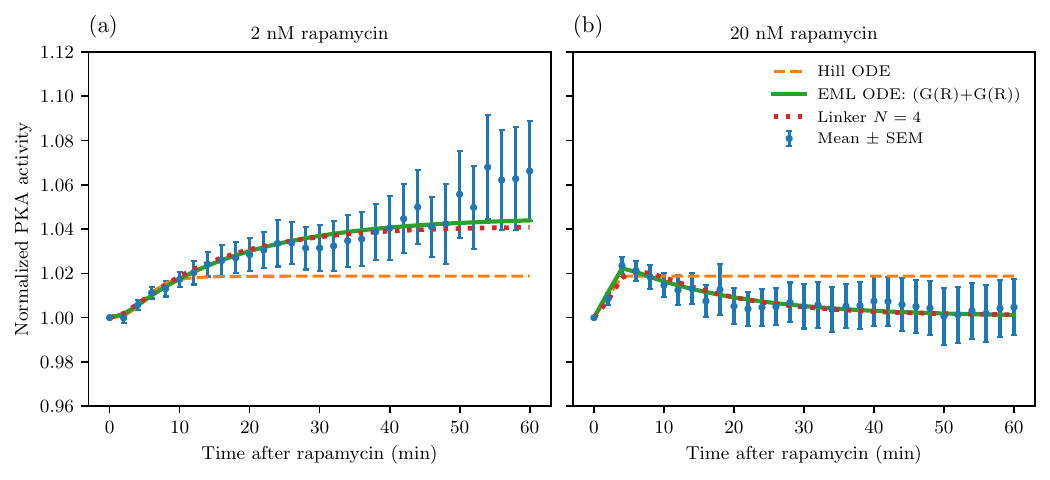}
		\caption{\label{fig:lacroix_fit}
		Restricted EML-ODE grammar search applied to PKA-R
		relocalization data from LaCroix \emph{et al.}~\cite{lacroix2022} Fig.~2F.
		(a)~$2$~nM rapamycin; (b)~$20$~nM rapamycin.
		Points show mean normalized plasma-membrane PKA activity $\pm$~SEM after rapamycin addition.
		The dashed curve is the Hill ODE; the dotted curve is the $N=4$ linker ODE from the LaCroix Appendix model; the
		solid curve is the best expression from an exhaustive search over all restricted EML-grammar ODEs with depth~$\le
		2$ and $\le 5$ nodes ($11$ expressions), which selected $E_*(R)=G_1(R)+G_2(R)$.
		All three models share the first-order relaxation structure $\tau\,\dot{y}_D = -y_D + 1 + B\,E(R_D(t))$,
		$y_D(0)=1$, with input $R_D(t)=D(1-e^{-kt})$ and centered gates $G_j(x)=(c_j+x)^{a_j}-b_j x - c_j^{a_j}$.
		The best fit parameters are recorded in Appendix~\ref{app:bestfit}.}
	\end{figure*}

	Figure~\ref{fig:lacroix_fit} shows the resulting two-dose fit. The Hill ODE remains structurally limited because it approaches a dose-dependent monotone plateau. It therefore misses the coexistence of a sustained low-dose response and a transient high-dose response. The linker model gives an excellent mechanistic description, as expected. The best EML-grammar ODE selected from the restricted search was
	\begin{equation}
		\label{eq:lacroix_selected}
		E_*(R)=G_1(R)+G_2(R),
	\end{equation}
	or, explicitly,
	\begin{equation}
		\label{eq:lacroix_selected_expanded}
		E_*(R)=
		(c_1+R)^{a_1}-b_1R-c_1^{a_1}
		+
		(c_2+R)^{a_2}-b_2R-c_2^{a_2}.
	\end{equation}
	Thus, the grammar selected a sum of two centered one-gate activation-suppression components sharing the same recruitment variable. The numerical ordering should not be overinterpreted: the EML expression has more fitted parameters than the linker model, and the linker model remains the mechanistically preferred description. The important point is that a restricted EML grammar is compatible with the known activation-suppression modeling and clearly improves on the Hill null.

\subsection{Grammar search on unresolved dynamics}
	\label{sec:nanda_grammar}

	The second example uses the Rho-GTPase recruitment experiment of Nanda \emph{et al.}, in which constitutively active Rho-family GTPases are acutely recruited to the plasma membrane and downstream activity sensors are followed over time \cite{nanda2023}. Unlike the previous example, these data are not accompanied by a compact model analogous to the LaCroix linker equation. They are therefore a better test of the methodological claim: EML as a grammar for proposing and selecting reduced equations.

\begin{figure*}[t]
		\centering
		\includegraphics[width=1.98\columnwidth]{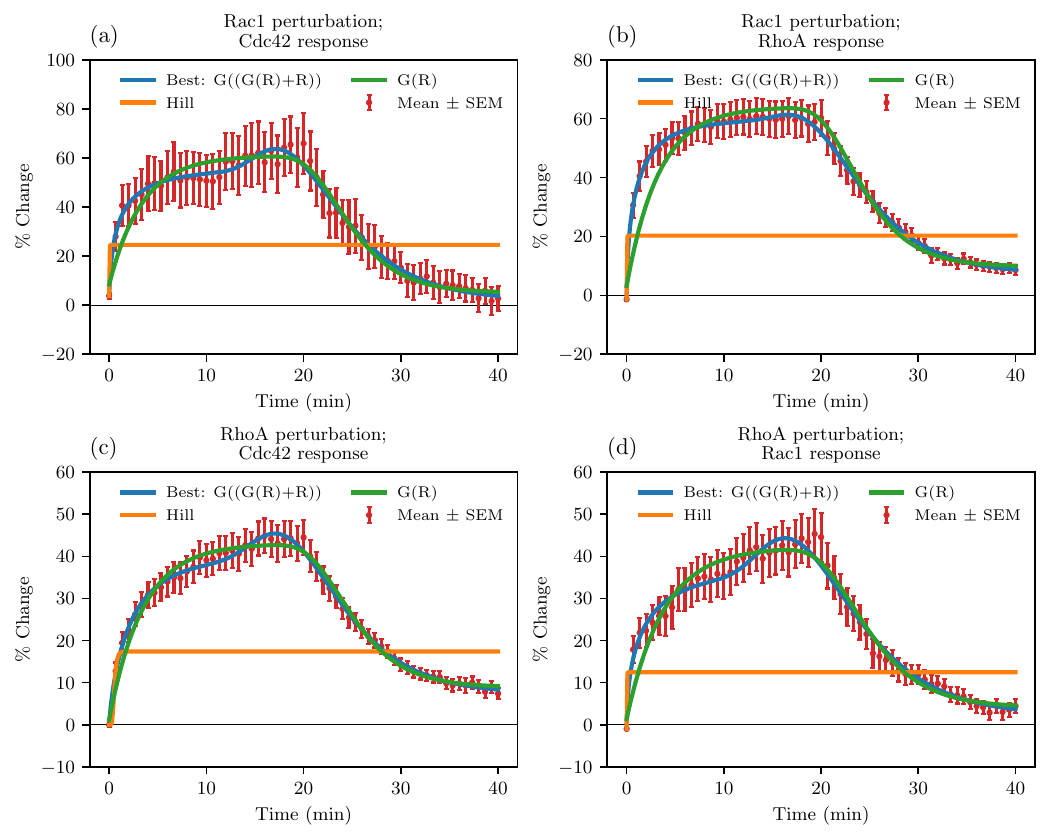}
		\caption{\label{fig:nanda_grammar}
		Restricted EML grammar search applied to all four response traces in
		figure~2d of Nanda \emph{et al.}~\cite{nanda2023}.
		Panels (a)-(d) show, respectively: Rac1 perturbation-Cdc42 response,
		Rac1 perturbation-RhoA response, RhoA perturbation-Cdc42 response, and
		RhoA perturbation-Rac1 response.
		Every second measured time point is plotted; error bars are SEM.
		Each panel shows the monotone Hill comparator, the one-gate expression
		$G(R)$, and the best expression $E_*$ found by exhaustive search over
		EML-grammar expressions with depth $\le 3$ and $\le 5$ nodes (see legend).
		The model is $y(t)=y_0+B\,E_*(R(t;k))$ with $R(t;k)=1-e^{-kt}$ and
		centered gates $G_j(x)=(c_j+x)^{a_j}-b_j\,x-c_j^{a_j}$. The best fit parameters are recorded in Appendix~\ref{app:bestfit}}
	\end{figure*}

	We applied the restricted grammar in Eq.~\eqref{eq:grammar} to the perturbations in RhoA and Rac1 and resulting response traces. The search was restricted to be small: expressions were enumerated up to EML depth three and five total expression-tree nodes, and continuous parameters were fitted on training time points. Candidate models were ranked by weighted error on held-out time points (Fig.~\ref{fig:nanda_grammar}). Two features of this search are worth noting. First, even at depth one the EML expression $G(R)$ (green curve in each panel) is already non-monotone and qualitatively tracks the rise-and-fall of the data, in contrast to the Hill comparator (orange) which saturates. A single Hill block cannot rise and then fall under a monotone recruitment input, whereas a single centered EML gate can. Second, the same depth-2 expression $G(G(R){+}R)$ is selected as the best held-out model for all four response traces. The algorithm did not merely compare a manually chosen Hill curve with a manually chosen EML curve; it searched a small EML expression space and the same branching expression won across four physically distinct perturbation-response combinations. We do not interpret this as a microscopic mechanism, but it does suggest that a common compositional motif: input combined with delayed-and-saturated input, then suppressed, captures the dominant temporal structure of these responses.

\begin{figure*}[t]
		\centering
		\includegraphics[width=1.98\columnwidth]{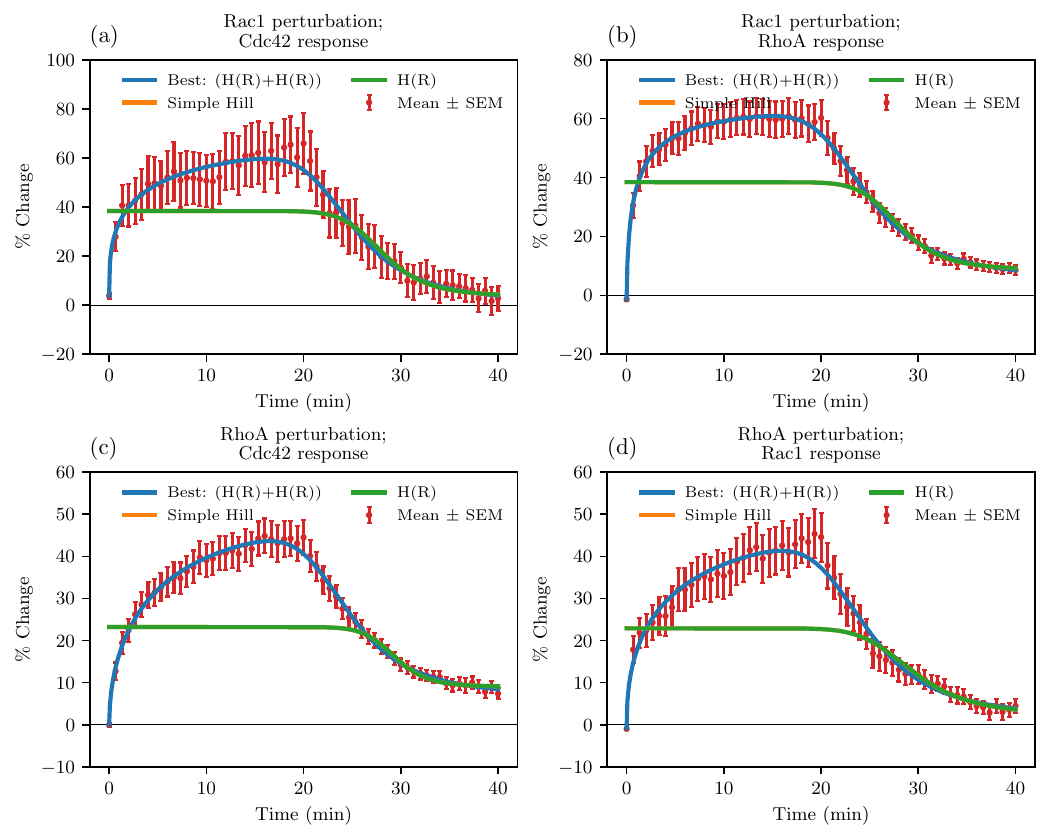}
		\caption{\label{fig:nanda_hill_grammar}
		Parallel Hill-grammar search applied to the same four response traces of Nanda \emph{et al.}~\cite{nanda2023} as in Fig.~\ref{fig:nanda_grammar}. The grammar is $E ::= R \mid H(E) \mid E+E$, with $H(\cdot)$ given by Eq.~\eqref{eq:hill_static}; the kinetic embedding, search bounds, and held-out validation procedure are identical to those used for the EML grammar in Fig.~\ref{fig:nanda_grammar}. ``Simple Hill'' is the depth-zero comparator with a single Hill block fit directly to the data, $H(R)$ is the best depth-one expression, and the curve labeled ``best'' is the best expression with depth $\le 3$ and $\le 5$ nodes. In every panel the best expression is $\bigl(H(R)+H(R)\bigr)$, a sum of two Hill blocks with opposing amplitudes, which closely matches the depth-2 EML fit in Fig.~\ref{fig:nanda_grammar}. At depth one, however, no Hill expression can produce the rise-then-fall: the orange ``Simple Hill'' and green $H(R)$ curves saturate to a monotone plateau. This contrasts with the depth-one EML expression $G(R)$ in Fig.~\ref{fig:nanda_grammar}, which is already non-monotone.}
	\end{figure*}

	To compare fairly with Hill-type modeling, we ran a parallel grammar search in which the EML gate was replaced by a Hill block of the same form, $E ::= R \mid H(E) \mid E+E$, with $H(\cdot)$ given by Eq.~\eqref{eq:hill_static}, and otherwise identical kinetic embedding, search bounds, and held-out validation. The results are shown in Fig.~\ref{fig:nanda_hill_grammar}. At depth one, $H(R)$ produces a monotone saturating curve in all four panels, as expected, and visibly fails to capture the late-time decay. The best Hill-grammar expression at depth $\le 3$ is $\bigl(H(R)+H(R)\bigr)$ in all four panels, and matches the EML fit in Fig.~\ref{fig:nanda_grammar} closely. The qualitative picture is therefore: at depth one, only the EML grammar produces a non-monotone response; at depth two, a sum of two Hills \cite{Vogel2016} (with opposing amplitudes) reproduces the EML fit. The EML grammar is not uniquely expressive, but it does collapse a depth-2 Hill construction into a depth-1 expression with fewer static-block parameters.

	The asymmetry can be quantified using AIC and BIC (Table~\ref{tab:nanda_aic_bic} in Appendix~\ref{app:bestfit}). At depth one, the EML gate $G(R)$ is preferred over the Hill block $H(R)$ in every panel by very large margins ($\Delta\text{AIC}$ between $277$ and $373$, $\Delta\text{BIC}$ between $275$ and $370$), reflecting that no single saturating block can fit the rise-then-fall shape. At depth two, the comparison narrows: the EML expression $G(G(R){+}R)$ ($p=9$) achieves the lowest AIC and BIC in panels~(a), (c), and~(d) despite carrying one more parameter than $\bigl(H(R)+H(R)\bigr)$ ($p=8$, with shared rate constant), with $\Delta\text{AIC}$ of $27$, $5.0$, and $10.0$ respectively. In panel~(b) the two grammars are essentially indistinguishable on held-out wMSE ($0.0698$ vs.\ $0.0737$), and BIC's parsimony preference selects the simpler Hill. Both information criteria therefore confirm what the curves in Fig.~\ref{fig:nanda_grammar} and Fig.~\ref{fig:nanda_hill_grammar} suggest visually: at depth two the two grammars are interchangeable in terms of fit quality, and the EML choice is justified primarily by parsimony at low depth and by the interpretability of each block as a single activation-suppression unit.

	\subsection{A toy coarse-graining benchmark}
	\label{sec:toy_coarse_graining}

	The grammar-search example above shows that EML expressions can be selected from data, but it does not demonstrate why higher EML depth could be useful. To test this point, we constructed a toy high-dimensional dynamical system whose output contains hidden delay structure. The microscopic model has $50$ internal states: a fast activation branch with $20$ first-order stages and a slower inhibitory branch with $30$ first-order stages. The imposed input is again a monotone recruitment variable $R(t)=1-\ee^{-k_R t}$. The measured output is a saturating positive function of the terminal activation state minus a saturating negative function of the terminal inhibitory state. Thus the overshoot and later adaptation arise from the separation of time scales in a high-dimensional network, not from an EML construction. The full ODE system used to generate the ground-truth trajectory is given in Appendix~\ref{app:toy_benchmark}.

\begin{figure*}[t]
		\centering
		\includegraphics[width=0.98\textwidth]{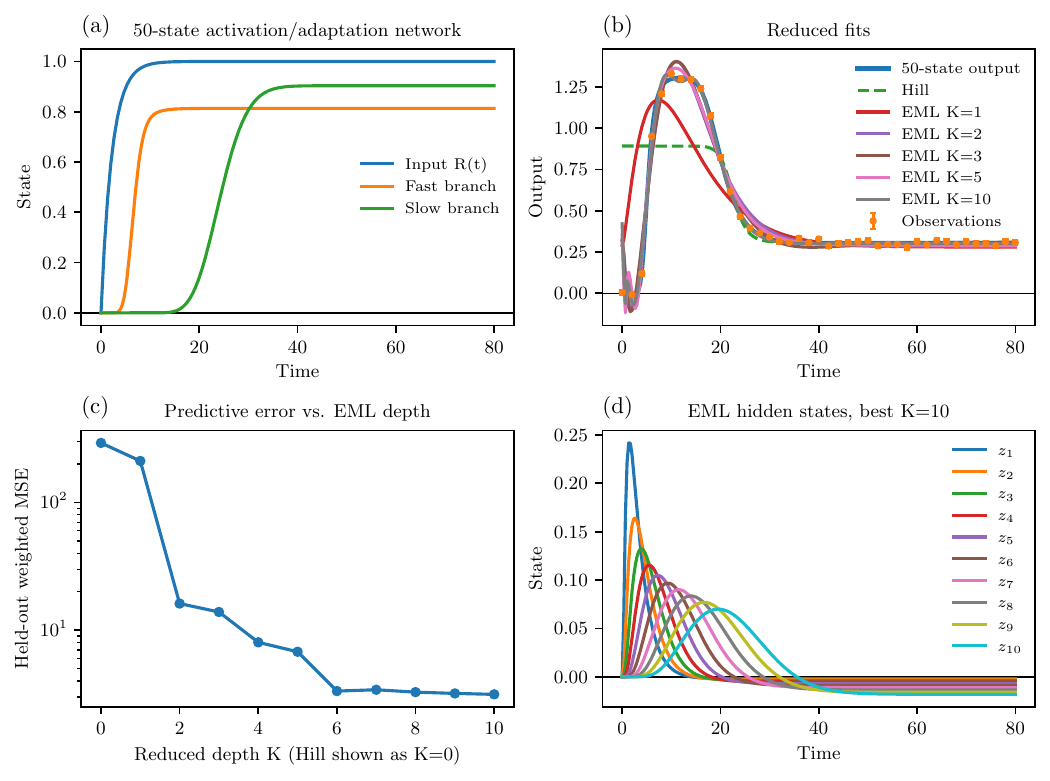}
		\caption{\label{fig:toy_coarse}
		Toy coarse-graining benchmark. A $50$-state activation/adaptation network driven by $R(t)=1-e^{-k_R t}$ is
		approximated by EML cascades of increasing depth $K$. Network parameters: $n_A=20$ fast states, $n_I=30$ slow
		states, $k_R=0.45$; fast branch $k_{\mathrm{on}}^A=2.4$, $k_{\mathrm{off}}^A=0.55$, $\tau_A=0.28$; slow branch
		$k_{\mathrm{on}}^I=0.75$, $k_{\mathrm{off}}^I=0.08$, $\tau_I=0.75$; output readout
		$y=1.6,A_{20}/(0.18+A_{20})-1.25,I_{30}/(0.22+I_{30})$; additive Gaussian noise $\sigma=0.015$ (seed 1). Each EML
		reduction generates $K$ centered hidden states from the same input and fits a linear readout by weighted least
		squares; every fourth time point is withheld for validation.
		(a) Trajectories of the full network: input $R(t)$, terminal fast-branch state $A_{20}(t)$, and terminal
		slow-branch state $I_{30}(t)$.
		(b) Reduced-model fits superimposed on the $50$-state output and noisy observations. The Hill model [$y_0=0.892$,
		$A=-0.585$, $K_d=0.327$, $h=20$, $k=0.018$; validation wMSE $=292$] and the $K=1$ EML cascade (wMSE $=211$)
		track the input monotonically and fail to reproduce the delayed adaptive peak. Increasing the EML depth
		progressively recovers the transient (wMSE $=16.2$ at $K=2$, $13.9$ at $K=3$, $6.79$ at $K=5$, and $3.15$ at
		$K=10$).
		(c) Held-out weighted MSE versus cascade depth $K$ (Hill model shown at $K=0$, log scale). The sharp drop between
		$K=1$ and $K=2$ marks the minimum depth required to represent the inhibitory delay.
		(d) Hidden states $z_1,\ldots,z_{10}$ of the best $K=10$ reduction form a sequence of progressively delayed,
		pulse-like basis functions.}
	\end{figure*}

	Can this $50$-state input-output map be approximated by a much lower-dimensional EML cascade? This benchmark was implemented in a reservoir-computing-like form \cite{maass2002,jaeger2004,nakajima2020}. In this paradigm, a complex nonlinear dynamical system (the 'reservoir') is kept fixed to generate a rich set of temporal basis functions, and only a simple linear readout layer is trained. By adopting this reservoir framework for our EML cascade, we treat the hidden states as a fixed basis. Thus, the calculation did not require optimizing all $a_k$, $b_k$, $c_k$, and $\tau_k$ parameters of a deep nonlinear cascade. For this benchmark the EML hidden states were generated by the centered cascade
	\begin{equation}
		\label{eq:toy_cascade}
		\tau_k \dot z_k=-z_k+
		G_{a_k,b_k,c_k}(z_{k-1}),
		\qquad z_0(t)=R(t),
	\end{equation}
	with $G_{a,b,c}$ defined in Eq.~\eqref{eq:centered_gate}. The readout was not chosen by hand. For each depth $K$, we fitted only a linear output layer,
	\begin{equation}
		\label{eq:toy_readout}
		y_K(t)=\beta_0+\sum_{j=1}^{K}\beta_j z_j(t),
	\end{equation}
	on training time points and evaluated the prediction on held-out time points. This makes the benchmark a conservative coarse-graining test: the EML cascade supplies a fixed structured temporal basis, while only the output weights are inferred from data. This choice reduces computational cost, avoids deep nonlinear convergence issues, and limits overfitting relative to fitting every cascade parameter. The results are shown in Fig.~\ref{fig:toy_coarse}. A monotone Hill response fails because the target output is generated by delayed activation followed by delayed inhibition. A one-state EML reduction is also insufficient. However, the held-out weighted mean squared error decreases sharply with EML depth. The largest improvement occurs between $K=1$ and $K=2$, indicating that one hidden EML state cannot represent the delayed inhibitory component, whereas a second layer already captures much of the activation-adaptation structure. Beyond $K=6$ the held-out wMSE plateaus, and AIC and BIC both select $K=6$ as the optimal depth (Table~\ref{tab:toy_aic_bic} in Appendix~\ref{app:toy_benchmark}); the additional layers $K=7,\ldots,10$ reduce held-out wMSE only marginally and are penalized by both information criteria.

	This benchmark shows that  EML depth can act as a coarse-grained dynamical basis, not that it recovers the microscopic species. Successive states $z_1,\ldots,z_K$ form delayed pulse-like components generated from the same monotone input, and the fitted readout combines these components to approximate the output of a much larger system. This is the operational meaning of higher-order EML in the present framework: it is a controlled way to increase reduced dynamical depth when hidden delays or distributed adaptive processes are present.

\subsection{Higher-order EML: solvable reduced ODE hierarchies}
	\label{sec:higher_order_eml}

	The LaCroix and Nanda examples use a first-order grammar search, whereas the toy benchmark uses a cascade of EML-generated states. The general construction is triangular, meaning each layer depends only on the output of the preceding layer, allowing the system to be solved sequentially. Let $z_0(t)=R(t)$ and define
	\begin{equation}
		\label{eq:z1_centered}
		\tau_1\dot z_1=-z_1+F_1(z_0),
	\end{equation}
	\begin{equation}
		\label{eq:zk_centered}
		\tau_k\dot z_k=-z_k+F_k(z_{k-1}),\qquad k=2,\ldots,K,
	\end{equation}
	where each drive is a centered EML gate
	\begin{equation}
		\label{eq:centered_drive}
		F_k(z)=(c_k+z)^{a_k}-c_k^{a_k}-b_kz
		=
		\eml(a_k\ln[c_k+z],\ee^{b_kz})-c_k^{a_k}.
	\end{equation}
	The subtraction of $c_k^{a_k}$ fixes the zero-input baseline, so that downstream layers are driven only by upstream activity. A measured response can then be modeled either as a terminal output, $y_K(t)=z_K(t)$, or as an observational readout
	\begin{equation}
		\label{eq:eml_linear_readout}
		y_K(t)=\beta_0+\sum_{j=1}^{K}\beta_jz_j(t),
	\end{equation}
	as in the toy coarse-graining benchmark.

	 Equations~(\ref{eq:z1_centered}--\ref{eq:centered_drive})  define a feedforward Hammerstein-Wiener cascade with first-order linear blocks; the contribution here is the use of the centered EML gate as the static nonlinearity. Compared to standard choices (sigmoid, Hill, polynomial), the EML gate provides non-monotone activation-suppression with a  closed-form derivative.	The hierarchy is exactly solvable:
	\begin{equation}
		\label{eq:general_convolution}
		z_k(t)=z_k(0)\ee^{-t/\tau_k}
		+\frac{1}{\tau_k}\int_0^t \ee^{-(t-s)/\tau_k}
		F_k(z_{k-1}(s))\,\dd s .
	\end{equation}
	Thus $z_1$ is solved from the known input, $z_2$ from $z_1$, and so on. General elementary closed forms are not expected for arbitrary $K$, but the convolution representation makes the reduced dynamics well defined and numerically stable.

	Linearization around a working point gives (Appendix~\ref{app:linear_response}),
	\begin{equation}
		\delta \dot z_k=-\frac{1}{\tau_k}\delta z_k
		+\frac{g_k}{\tau_k}\delta z_{k-1},
		\qquad
		g_k=a_k(c_k+z_{k-1}^*)^{a_k-1}-b_k .
	\end{equation}
	The corresponding transfer function is
	\begin{equation}
		\label{eq:transfer}
		H_K(s)=\prod_{k=1}^{K}\frac{g_k}{1+s\tau_k}.
	\end{equation}
	Each additional layer therefore contributes one gain factor and one timescale. In this reduced sense, EML depth is a model-selection parameter for hidden delay and adaptive structure, not a direct count of molecular regulatory steps.

\section{Discussion}
	
	The examples above suggest a narrow but useful role for EML-generated ODEs as a model-discovery grammar. The framework is best summarized by three claims.

	First, the centered EML gate of Eq.~\eqref{eq:centered_gate} is a minimal non-monotone activation-suppression block, whereas any block built from monotone saturating primitives such as the Hill function requires at least a difference of two such blocks with opposing amplitudes to produce the same shape. The Nanda \emph{et al.} analysis (Fig.~\ref{fig:nanda_grammar}, Fig.~\ref{fig:nanda_hill_grammar}) makes this concrete: at depth one.

	Second, at sufficient compositional depth other grammars match EML. A Hill grammar at depth two reproduces the EML fit on all four Nanda response traces using a sum of two Hills with opposing amplitudes. The quantitative comparison in Sec.~\ref{sec:nanda_grammar} (AIC/BIC, Table~\ref{tab:nanda_aic_bic}) shows that the depth-1 EML gate $G(R)$ beats the depth-1 Hill block $H(R)$ by very large margins in all four panels, while at depth two the two grammars become essentially interchangeable, with EML preferred by AIC/BIC in three panels and Hill preferred in the fourth. The EML grammar is therefore not uniquely expressive in the limit of unrestricted depth; its advantage is parsimony at low depth and interpretability of the individual block as a single activation-suppression unit.

	Third, the centered EML cascade of Eqs.~(\ref{eq:z1_centered}--\ref{eq:centered_drive}) provides a controlled ladder between simple Hill-type phenomenology and detailed mechanistic models. The cascade is a feedforward Hammerstein-Wiener architecture in which model depth corresponds to hidden activation-suppression regulatory depth; each additional layer contributes one gain factor and one timescale to the linearized transfer function (Eq.~\eqref{eq:transfer}). On the toy benchmark this ladder compresses a $50$-state activation-adaptation network into a low-dimensional reservoir with a learned linear readout. AIC and BIC both select $K=6$ as the optimal cascade depth (Table~\ref{tab:toy_aic_bic}), one depth below the held-out wMSE plateau, identifying an effective coarse-grained dimensionality of six for this network.

	When an established mechanism exists, as in the PKA-R linker model \cite{lacroix2022}, the mechanistic model is preferable. The EML model is valuable there because it demonstrates that a one-gate activation-suppression module recovers a real biochemical response motif at held-out wMSE comparable to the linker model, even though the linker model is the mechanistically preferred description.
	
	Several limitations are immediate. First, as the Hill-grammar comparison shows, the EML hierarchy is not uniquely expressive; other grammars can generate biphasic and adaptive responses at modestly higher compositional depth. Second, the biological interpretation of $a_k$, $b_k$, and $c_k$ is reduced-model interpretation, not direct molecular measurement unless additional experiments constrain the underlying processes. Third, increasing $K$ can overfit unless penalized by predictive validation or information criteria such as AIC or BIC (Appendix~\ref{app:computational_grammar}). Selection in this manuscript uses held-out wMSE as the primary criterion; AIC and BIC are reported alongside (Tables~\ref{tab:nanda_aic_bic} and~\ref{tab:toy_aic_bic}) and broadly confirm the held-out wMSE rankings, with the one exception of Nanda panel~(b) noted above. These limitations are not defects specific to EML; they are the usual constraints on reduced dynamical modeling.
	
	The proposed value of EML is that it organizes model expansion and symbolic regression. Instead of adding arbitrary empirical terms until a curve is fit, one increases EML depth and asks whether the data justify another activation-suppression layer. If the answer is no, the first-order model is sufficient. If the answer is yes, the added layer has a clear dynamical interpretation and an exact convolution solution. This makes EML-generated ODEs a plausible grammar for model discovery in systems where perturbations induce overshoot, adaptation, or biphasic responses but a detailed mechanistic model is unavailable.

	The universality of EML is a risk as well as a strength. Without a restricted grammar and explicit penalties, EML can generate arbitrary curve-fitting expressions. The framework is therefore meaningful only when the allowed terminals, EML depth, parameter count, and validation procedure are specified in advance. In an era
	increasingly dominated by opaque machine learning mod-
	els, EML offers a mathematically transparent bridge be-
	tween data-driven discovery and interpretable nonlinear
	dynamics.

	\section{Acknowledgments}
	Numerical simulations were paid for by AE's startup funds.
	
    \subsection*{Data availability}
    All scripts, processed inputs, and commands required to reproduce the figures are available at \url{https://www.github.com/AmirErez/Manuscript_EML_biophysics}.

	%\bibliography{eml}
	
	\bibliographystyle{unsrt}

	\clearpage
	\appendix
	\part*{Appendix}
	
	\pagenumbering{arabic}
	\renewcommand{\thefigure}{S\arabic{figure}}
	\setcounter{figure}{0} 
	\renewcommand{\thetable}{S\arabic{table}}
	\renewcommand{\thesubsection}{S\arabic{subsection}}
	\renewcommand{\theequation}{S\arabic{equation}}
	\setcounter{equation}{0} 
	
	\section{Derivation of the one-gate EML activation-suppression module}
	\label{app:module}
	
	Starting from
	\begin{align}
		\eml(\alpha\ln R,\ee^{\beta R})
		&=\ee^{\alpha\ln R}-\ln\left(\ee^{\beta R}\right) \\
		&=R^\alpha-\beta R .
	\end{align}
	For $R>0$, the derivative is
	\begin{equation}
		\frac{\dd}{\dd R}\left(R^\alpha-\beta R\right)
		=\alpha R^{\alpha-1}-\beta .
	\end{equation}
	For $0<\alpha<1$ and $\beta>0$, the derivative vanishes at
	\begin{equation}
		R_*^{\alpha-1}=\frac{\beta}{\alpha},
	\end{equation}
	so
	\begin{equation}
		R_*=\left(\frac{\alpha}{\beta}\right)^{1/(1-\alpha)} .
	\end{equation}
	The second derivative is
	\begin{equation}
		\frac{\dd^2}{\dd R^2}\left(R^\alpha-\beta R\right)
		=\alpha(\alpha-1)R^{\alpha-2}<0,
	\end{equation}
	so $R_*$ is a maximum.
	
	\section{Solution of the first-order EML relaxation equation}
	\label{app:solutions}
	
	Consider
	\begin{equation}
		\tau\dot y=-y+F(t),
	\end{equation}
	where
	\begin{equation}
		F(t)=y_0+B\left[R(t)^\alpha-\beta R(t)\right].
	\end{equation}
	The ODE can be written as
	\begin{equation}
		\dot y+\frac{1}{\tau}y=\frac{1}{\tau}F(t).
	\end{equation}
	Multiplying by the integrating factor $\ee^{t/\tau}$ gives
	\begin{equation}
		y(t)=y(0)\ee^{-t/\tau}
		+\frac{1}{\tau}\int_0^t \ee^{-(t-s)/\tau}F(s)\,\dd s .
	\end{equation}
	
	\section{Best fit parameters for the three figures}
		\label{app:bestfit}

		Tables~\ref{tab:lacroix_bestfit} and~\ref{tab:nanda_bestfit} report the best-fit parameters
		for the EML-grammar fits in Figs.~\ref{fig:lacroix_fit} and~\ref{fig:nanda_grammar},
		respectively. Held-out wMSE values for the comparator models are listed alongside.

		\begin{table*}[h!]
		\centering
		\small
		\caption{Best-fit parameters for Fig.~\ref{fig:lacroix_fit} (LaCroix \emph{et al.}\
		PKA-R relocalization). Selected EML-ODE expression: $E_*=G_1(R)+G_2(R)$
		(depth~1, 5~nodes, 2~gates, 9~parameters). Comparator models share the same
		first-order relaxation structure. wMSE is reported as train/hold.}
		\label{tab:lacroix_bestfit}
		\begin{tabular}{lll}
		\hline
		Model & Parameters & wMSE \\
		\hline
		EML ODE   & $B=10$, $k=1.365$~min$^{-1}$, $\tau=16.6$~min;        & $0.131$ / $0.201$ \\
		$G_1(R)+G_2(R)$ & gate~1: $(a_1,b_1,c_1)=(1.02,\,0.517,\,0.0611)$;     &  \\
		          & gate~2: $(a_2,b_2,c_2)=(0.808,\,1.005,\,10.5)$        &  \\
		\hline
		Linker $N=4$ & $A=6.95{\times}10^5$, $S_0=0.833$,                   & $0.208$ / $0.210$ \\
		             & $q=1.35{\times}10^{-4}$,                              &  \\
		             & $k=0.503$~min$^{-1}$, $\tau=12.95$~min                 &  \\
		\hline
		Linker $N=2$ & $A=2.71{\times}10^5$, $S_0=1.41$,                    & $0.208$ / $0.210$ \\
		             & $q=2.96{\times}10^{-4}$,                              &  \\
		             & $k=0.503$~min$^{-1}$, $\tau=12.94$~min                 &  \\
		\hline
		Hill ODE     & $A=0.0188$, $K_d=0.00851$, $h=3.45$,                  & $1.70$  / $1.62$  \\
		             & $k=0.001$~min$^{-1}$, $\tau=0.081$~min                 &  \\
		\hline
		\end{tabular}
		\end{table*}

		\begin{table*}[h!]
		\centering
		\small
		\caption{Best-fit parameters for the EML-grammar selection in Fig.~\ref{fig:nanda_grammar}
		(Nanda~\emph{et~al.}\ Rho-GTPase recruitment). All four panels select
		$E_*=G(G(R){+}R)$ (depth~2, 5~nodes, 2~gates, 9~parameters).
		Held-out (validation) wMSE for the selected expression, for the depth-1 EML
		expression $G(R)$, and for the Hill comparator are listed for reference.}
		\label{tab:nanda_bestfit}
		\begin{tabular}{lllllllll}
		\hline
		Panel & $y_0$ & $B$ & $k$~(min$^{-1}$) & Outer gate $(a_1,b_1,c_1)$ & Inner gate $(a_2,b_2,c_2)$
		      & wMSE [$G(G(R){+}R)$] & wMSE [$G(R)$] & wMSE [Hill] \\
		\hline
		(a) & $3.82$   & $165$  & $0.205$ & $(0.110,\,0.0378,\,0.0207)$       & $(71.2,\,4.5{\times}10^{-6},\,0.042)$  & $0.128$  & $0.236$ & $6.38$ \\
		(b) & $-1.37$  & $221$  & $0.213$ & $(0.506,\,0.868,\,1.8{\times}10^{-8})$ & $(74.9,\,0.765,\,3.3{\times}10^{-5})$  & $0.0698$ & $1.16$  & $28.0$ \\
		(c) & $0.005$  & $745$  & $0.211$ & $(0.966,\,0.947,\,{\approx}0)$    & $(54.5,\,1.1{\times}10^{-6},\,0.025)$  & $0.161$  & $0.474$ & $23.5$ \\
		(d) & $-0.949$ & $86.2$ & $0.188$ & $(0.481,\,0.475,\,6.6{\times}10^{-10})$ & $(36.1,\,0.664,\,0.037)$              & $0.234$  & $0.736$ & $14.4$ \\
		\hline
		\end{tabular}
		\end{table*}

		\subsubsection*{AIC and BIC diagnostics (Fig.~\ref{fig:nanda_grammar})}
		\label{app:nanda_aic_bic}
		All four panels share $N_T=91$ training points (alternating split from 121 total time points).
		Table~\ref{tab:nanda_aic_bic} reports AIC and BIC for the Hill comparator $H(R)$ ($p=5$),
		the depth-1 EML gate $G(R)$ ($p=6$), the double-Hill comparator $(H(R){+}H(R))$ ($p=8$,
		shared~$k$), and the best EML expression $G(G(R){+}R)$ ($p=9$), all evaluated on identical
		training residuals using Eq.~\eqref{eq:aic_bic}. $\Delta$AIC and $\Delta$BIC are relative
		to the per-panel AIC minimum.
		In panels~(a), (c), and~(d) the EML expression $G(G(R){+}R)$ achieves the lowest AIC and BIC.
		In panel~(b) the EML and Hill grammars are essentially indistinguishable on held-out
		wMSE ($0.0698$ vs.\ $0.0737$, a difference within sampling noise), and BIC's parsimony
		preference selects the simpler double-Hill expression
		($\Delta$AIC $=4.3$ for $G(G(R){+}R)$). In all panels $G(R)$ ($p=6$) is substantially
		preferred over $H(R)$ ($p=5$) on both AIC and BIC, confirming the parsimony advantage of
		the single EML gate. The double-Hill fit in panel~(c) has $h_2=30$ at its upper search
		bound, so its AIC should be interpreted with caution.

		\begin{table*}[h!]
		\centering
		\small
		\caption{AIC and BIC for selected models fitted to the four Nanda~\emph{et~al.}\ response
		traces (Fig.~\ref{fig:nanda_grammar}). $N_T=91$; $\chi^2_T=N_T\times\mathrm{wMSE}_{\rm train}$.
		$\Delta$AIC and $\Delta$BIC are relative to the per-panel AIC minimum (marked with~$\star$).
		$\dagger$: $h_2=30$ at upper search bound.}
		\label{tab:nanda_aic_bic}
		\begin{tabular}{llrrrrrrr}
		\hline
		Panel & Model & $p$ & $\chi^2_T$ & wMSE$_{\rm hold}$ & AIC & BIC & $\Delta$AIC & $\Delta$BIC \\
		\hline
		(a) & $G(G(R){+}R)$$^\star$ & 9 &  6.9 & 0.128  & $-216$ & $-194$ &   0  &   0 \\
		(a) & $(H(R){+}H(R))$       & 8 &  9.5 & 0.260  & $-189$ & $-169$ &  27  &  24 \\
		(a) & $G(R)$                & 6 & 26.7 & 0.236  & $-100$ &  $-85$ & 116  & 109 \\
		(a) & $H(R)$                & 5 & 571  & 6.38   &   177  &   190  & 393  & 384 \\
		\hline
		(b) & $(H(R){+}H(R))$$^\star$ & 8 &  6.5 & 0.0737 & $-224$ & $-204$ &  0.0 &  0.0 \\
		(b) & $G(G(R){+}R)$          & 9 &  6.7 & 0.0698 & $-219$ & $-197$ &  4.3 &  6.8 \\
		(b) & $G(R)$                  & 6 & 104  & 1.16   &    24  &    39  & 248  & 243  \\
		(b) & $H(R)$                  & 5 & 2580 & 28.0   &   314  &   327  & 538  & 530  \\
		\hline
		(c) & $G(G(R){+}R)$$^\star$ & 9 &  8.3 & 0.161  & $-200$ & $-177$ &  0.0 &  0.0 \\
		(c) & $(H(R){+}H(R))$$^\dagger$ & 8 &  9.0 & 0.0774 & $-195$ & $-175$ &  5.0 &  2.5 \\
		(c) & $G(R)$                & 6 & 34.7 & 0.474  &  $-76$ &  $-61$ & 124  & 116  \\
		(c) & $H(R)$                & 5 & 2127 & 23.5   &   297  &   309  & 496  & 487  \\
		\hline
		(d) & $G(G(R){+}R)$$^\star$ & 9 & 19.2 & 0.234  & $-124$ & $-101$ &  0.0 &  0.0 \\
		(d) & $(H(R){+}H(R))$       & 8 & 21.9 & 0.270  & $-114$ &  $-93$ & 10.0 &  7.5 \\
		(d) & $G(R)$                & 6 &  52  & 0.736  &  $-38$ &  $-23$ &  85  &  78  \\
		(d) & $H(R)$                & 5 & 1356 & 14.4   &   256  &   268  & 379  & 369  \\
		\hline
		\end{tabular}
		\end{table*}

	\section{Dynamical form of the LaCroix linker model}
	\label{app:linker}
	
	LaCroix \emph{et al.} model signaling as proportional to the concentration of fully occupied linker molecules \cite{lacroix2022}. Their heuristic expression can be written, up to scale, as a function of an effective linker variable $S$. For $N$ binding sites, the reduced form used here is
	\begin{equation}
		\Phi_N(S)=\frac{S}{S^N+\sum_{n=1}^{N} n S^{N-n}} .
	\end{equation}
	For $N=4$,
	\begin{equation}
		\Phi_4(S)=\frac{S}{S^4+S^3+2S^2+3S+4} .
	\end{equation}
	To construct a dynamical model, we assume that recruitment increases the local linker abundance as
	\begin{equation}
		S_D(t)=S_0+qD(1-\ee^{-k_Rt}).
	\end{equation}
	The observed PKA activity is then assumed to relax toward the shifted fully occupied linker level,
	\begin{equation}
		\tau\dot y_D=-y_D+1+A\left[\Phi_N(S_D(t))-\Phi_N(S_0)\right].
	\end{equation}
	The subtraction of $\Phi_N(S_0)$ enforces the normalized initial condition $y_D(0)=1$ when the data are normalized to their pre-rapamycin baseline.
	
	\section{Linear response of the centered higher-order EML cascade}
	\label{app:linear_response}
	Let $z_k^*$ be a fixed working point and define $\delta z_k=z_k-z_k^*$. Linearizing
	\begin{equation}
		\tau_k\dot z_k=-z_k+F_k(z_{k-1})
	\end{equation}
	gives
	\begin{equation}
		\tau_k\delta\dot z_k=-\delta z_k+F_k'(z_{k-1}^*)\delta z_{k-1} .
	\end{equation}
	Thus
	\begin{equation}
		\delta\dot z_k=-\frac{1}{\tau_k}\delta z_k
		+\frac{g_k}{\tau_k}\delta z_{k-1},
	\end{equation}
	where
	\begin{equation}
		g_k=F_k'(z_{k-1}^*).
	\end{equation}
	For $k\geq 2$,
	\begin{equation}
		F_k'(z)=a_k(c_k+z)^{a_k-1}-b_k.
	\end{equation}
	The centering constant $-c_k^{a_k}$ does not affect the derivative, but it changes the baseline and is therefore essential for constructing pulse-propagating cascades. Taking the Laplace transform gives
	\begin{equation}
		(s+1/\tau_k)\delta Z_k(s)=\frac{g_k}{\tau_k}\delta Z_{k-1}(s),
	\end{equation}
	or
	\begin{equation}
		\frac{\delta Z_k(s)}{\delta Z_{k-1}(s)}=\frac{g_k}{1+s\tau_k} .
	\end{equation}
	Multiplying over $k=1,\ldots,K$ gives
	\begin{equation}
		H_K(s)=\prod_{k=1}^{K}\frac{g_k}{1+s\tau_k} .
	\end{equation}
	This proves Eq.~\eqref{eq:transfer}.

	\section{Computational EML grammar and validation objective}
	\label{app:computational_grammar}
	For the Nanda analysis, all expressions with depth at most three and at most five tree nodes were enumerated. For the LaCroix analysis, the same grammar was used but embedded in a first-order relaxation equation so that the selected expression defined the target of an ODE rather than an instantaneous static response. The search space was deliberately small. Its purpose was to demonstrate the grammar and validation procedure, not to claim global optimality over all possible EML expressions.

	For a candidate expression $E$ with continuous parameters $\theta$, predictions $\hat y_E(t_i;\theta)$, observations $y_i$, and reported SEM values $\sigma_i$, the SEM values were floored as
	\begin{equation}
		\tilde\sigma_i=\max(\sigma_i,\sigma_{\rm floor}),
		\qquad
		\sigma_{\rm floor}=0.25\,{\rm median}\{\sigma_i:\sigma_i>0\}.
	\end{equation}
	Time points were split deterministically into training and held-out sets. Parameters were fitted on the training set,
	\begin{equation}
		\hat\theta_E=\arg\min_{\theta}
		\chi^2_{\mathcal T}(E,\theta),
		\qquad
		\chi^2_{\mathcal T}(E,\theta)=
		\sum_{i\in\mathcal T}
		\left[
		\frac{y_i-\hat y_E(t_i;\theta)}{\tilde\sigma_i}
		\right]^2 .
	\end{equation}
	Model selection used the held-out score
	\begin{equation}
		\mathcal S(E)=
		\frac{1}{|\mathcal V|}
		\sum_{i\in\mathcal V}
		\left[
		\frac{y_i-\hat y_E(t_i;\hat\theta_E)}{\tilde\sigma_i}
		\right]^2
		+\lambda_d d(E)+\lambda_n n(E),
	\end{equation}
	where $d(E)$ is EML depth and $n(E)$ is node count. The figures reported in the main text used $\lambda_d=\lambda_n=0$; the penalty terms are included to define the general search objective and can be activated in larger searches. AIC and BIC were also computed from the training residuals,
	\begin{equation}\label{eq:aic_bic}
		{\rm AIC}=N_{\mathcal T}\ln\left(\frac{\chi^2_{\mathcal T}}{N_{\mathcal T}}\right)+2p,
		\qquad
		{\rm BIC}=N_{\mathcal T}\ln\left(\frac{\chi^2_{\mathcal T}}{N_{\mathcal T}}\right)+p\ln N_{\mathcal T},
	\end{equation}
	where $p$ is the number of fitted parameters. These information criteria are diagnostic summaries rather than the primary selection rule in the figures shown.

\section{Toy coarse-graining benchmark}
	\label{app:toy_benchmark}

	The toy benchmark used a high-dimensional network that was not generated from EML. The imposed input was
	\begin{equation}
		R(t)=1-\exp(-k_Rt).
	\end{equation}
	The fast branch contained $n_A=20$ first-order states, denoted $A_1,\ldots,A_{n_A}$, and the slow branch contained $n_I=30$ first-order states, denoted $I_1,\ldots,I_{n_I}$. The first state in each branch obeyed an input-dependent activation equation,
	\begin{align}
		\dot A_1 &= k_{\rm on}^A R(t)(1-A_1)-k_{\rm off}^A A_1,\\
		\dot I_1 &= k_{\rm on}^I R(t)(1-I_1)-k_{\rm off}^I I_1.
	\end{align}
	Subsequent states formed two linear delay chains,
	\begin{align}
		\dot A_j &= \frac{A_{j-1}-A_j}{\tau_A},
		\qquad j=2,\ldots,n_A,\\
		\dot I_j &= \frac{I_{j-1}-I_j}{\tau_I},
		\qquad j=2,\ldots,n_I.
	\end{align}
	The ground-truth output was
	\begin{equation}
		y_{\rm true}(t)=
		y_0+
		A_{\rm amp}\frac{A_{n_A}(t)}{K_A+A_{n_A}(t)}
		-
		I_{\rm amp}\frac{I_{n_I}(t)}{K_I+I_{n_I}(t)}.
	\end{equation}
	Independent Gaussian noise was added to obtain the observed benchmark trace,
	\begin{equation}
		y_{\rm obs}(t_i)=y_{\rm true}(t_i)+\eta_i,
		\qquad
		\eta_i\sim\mathcal N(0,\sigma_{\rm noise}^2).
	\end{equation}
	The representative simulation used $n_A=20$, $n_I=30$, $k_R=0.45$, $k_{\rm on}^A=2.4$, $k_{\rm off}^A=0.55$, $\tau_A=0.28$, $k_{\rm on}^I=0.75$, $k_{\rm off}^I=0.08$, $\tau_I=0.75$, $A_{\rm amp}=1.6$, $I_{\rm amp}=1.25$, $K_A=0.18$, $K_I=0.22$, and $\sigma_{\rm noise}=0.015$.

	For the reduced model, a fixed centered EML cascade generated hidden states $z_1,\ldots,z_K$ from the same input:
	\begin{equation}
		\tau_k\dot z_k=-z_k+G_{a_k,b_k,c_k}(z_{k-1}),
		\qquad z_0(t)=R(t).
	\end{equation}
	The reservoir was deterministic rather than random. For layer $k=1,\ldots,K$, the fixed gate parameters were
	\begin{align}
		a_k &= 0.45+0.035(k-1),\\
		b_k &=
		\begin{cases}
			1.00, & k=1,\\
			0.42-0.015\min(k-1,10), & k\ge 2,
		\end{cases}\\
		c_k &=
		\begin{cases}
			10^{-6}, & k=1,\\
			0.08, & k\ge 2,
		\end{cases}
	\end{align}
	and the layer time constants were
	\begin{equation}
		\tau_k=\tau_0\left[1+0.55(k-1)\right].
	\end{equation}
	The input rate used by the reservoir was written as $R(t;k_{\rm fit})=1-\exp(-k_{\rm fit}t)$. The two global reservoir hyperparameters, $k_{\rm fit}$ and $\tau_0$, were selected by an explicit grid search, with $k_{\rm fit}$ taking $18$ equally spaced values between $0.15$ and $0.80$ and $\tau_0$ taking $20$ equally spaced values between $0.5$ and $5.5$. For each grid point and each depth $K$, only the readout coefficients in Eq.~\eqref{eq:toy_readout} were fitted by weighted least squares on training time points, and performance was measured on held-out time points using Eq.~\eqref{eq:validation_score}. Thus the benchmark evaluates whether an EML-generated temporal basis can approximate the input-output behavior of a much larger hidden network.

	\subsubsection*{AIC and BIC by cascade depth}
	\label{app:toy_aic_bic}
	For the EML cascade at depth $K$, the number of fitted parameters is $p=K+3$:
	$K+1$ readout coefficients $\beta_0,\ldots,\beta_K$ plus the two grid-searched reservoir
	hyperparameters $k_{\rm fit}$ and $\tau_0$. AIC and BIC are computed from training residuals
	using Eq.~\eqref{eq:aic_bic} with the same floored weights as the validation score; the Hill
	comparator uses $p=5$. Both criteria are minimized at $K=6$, one depth below the held-out wMSE
	plateau (which is nearly flat from $K=6$ to $K=10$, Table~\ref{tab:toy_aic_bic}).
	The abrupt drop in training $\chi^2$ between $K=1$ and $K=2$ (from $38\,124$ to $2\,765$)
	reflects the onset of the hidden slow branch: a single EML state cannot reproduce the
	adaptation, but two states can.

	\begin{table}[h!]
	\centering
	\small
	\caption{AIC and BIC for the toy cascade as a function of EML depth $K$
	(Fig.~\ref{fig:toy_coarse}c). $K=0$ is the Hill comparator ($p=5$); $K\ge1$ are
	EML cascades with $p=K+3$. $N_T=181$ for all rows.
	$\Delta$AIC and $\Delta$BIC are relative to the AIC minimum at $K=6$ (marked~$\star$).}
	\label{tab:toy_aic_bic}
	\begin{tabular}{rrrrrrrrrr}
	\hline
	$K$ & $p$ & $\chi^2_T$ & wMSE$_{\rm hold}$ & AIC & BIC & $\Delta$AIC & $\Delta$BIC \\
	\hline
	 0 &  5 & 51\,367 & 292   & 1032 & 1048 & 665 & 653 \\
	 1 &  4 & 38\,124 & 211   &  976 &  989 & 610 & 594 \\
	 2 &  5 &  2\,765 &  16.2 &  503 &  519 & 137 & 124 \\
	 3 &  6 &  2\,277 &  13.9 &  470 &  490 & 103 &  94 \\
	 4 &  7 &  1\,955 &   8.05 &  445 &  467 &  78 &  72 \\
	 5 &  8 &  1\,839 &   6.79 &  436 &  461 &  69 &  66 \\
	 6$^\star$ &  9 &  1\,244 &   3.34 &  367 &  396 &   0 &   0 \\
	 7 & 10 &  1\,239 &   3.42 &  368 &  400 & 1.3 & 4.5 \\
	 8 & 11 &  1\,288 &   3.28 &  377 &  412 &  10 &  17 \\
	 9 & 12 &  1\,282 &   3.21 &  378 &  417 &  11 &  21 \\
	10 & 13 &  1\,281 &   3.15 &  380 &  422 &  13 &  26 \\
	\hline
	\end{tabular}
	\end{table}
\end{document}